\documentclass[11pt,oneside,leqno]{amsart}

\usepackage{amsxtra}
\usepackage{amsopn}
\usepackage{amsmath,amsthm,amssymb}
\usepackage{color}
\usepackage{amscd}
\usepackage{pifont}
\usepackage{amsfonts}
\usepackage{latexsym}
\usepackage{verbatim}
\usepackage{pb-diagram}

\theoremstyle{plain}
\newtheorem{theorem}{Theorem}[section]
\newtheorem*{theorem*}{Theorem}

\newtheorem{cor}[theorem]{Corollary}

\newtheorem*{mt*}{Main Theorem}
\newcommand\g{{\mathfrak{g}}}

\begin{document}
\title[$3D$ homogeneous generalized Ricci solitons]{Three-dimensional homogeneous \\ generalized Ricci solitons}
\author{Giovanni Calvaruso}
\date{}

\address{Dipartimento di Matematica e Fisica \lq\lq E. De Giorgi\rq\rq \\
Universit\`a del Salento\\
Prov. Lecce-Arnesano \\
73100 Lecce\\ Italy.}
\email{giovanni.calvaruso@unisalento.it}

\thanks{Author partially supported by funds of the University of Salento and MIUR (PRIN)}

\begin{abstract}
We study three-dimensional generalized Ricci solitons, both in Riemannian and Lorentzian settings. We shall determine their homogeneous models, classifying left-invariant generalized Ricci solitons on three-dimensional Lie groups.
 \end{abstract}

\maketitle

\section{Introduction}

Generalized Ricci solitons were recently introduced in \cite{NR}.  A {\em generalized Ricci soliton} is a pseudo-Riemannian manifold $(M,g)$ admitting  a smooth vector field $X$, such that
\begin{equation}\label{GRS}
\mathcal{L}_X g+2 \alpha X^{\flat}\odot  X^\flat - 2\beta Ric = 2 \lambda g,
\end{equation}
for some real constants $\alpha,\beta,\lambda$, where $\mathcal{L}_X$ denotes the Lie derivative in the direction of $X$, $X^\flat$ denotes a $1$-form such that $X^\flat (Y)= g(X, Y)$ and $Ric$ is the Ricci tensor. 

For particular values of the constants $\alpha,\beta,\lambda$, several important equations occur as special cases of equation~\eqref{GRS}. In particular, one has:
\begin{itemize}
\item[(K)] the {\em Killing vector field equation} when $\alpha=\beta=\lambda=0$;
\item[(H)] the {\em homothetic vector field equation} when $\alpha=\beta=0$;
\item[(RS)] the {\em Ricci soliton equation} when $\alpha=0$ and $\beta =1$  \cite{Cao}.
\item[(E-W)] a special case of the {\em Einstein-Weyl equation} in conformal geometry when $\alpha=1$ and $\beta=-\frac{1}{n-2}$ ($n>2$) \cite{CP};
\item[(PS)] the equation for a {\em metric projective structure} with a {\em skew-symmetric Ricci tensor} representative in the projective class when  $\alpha=1$, $\beta=-\frac{1}{n-1}$ and  $\lambda =0$ \cite{Ra};
\item[(VN-H)] the {\em vacuum near-horizon geometry equation} of a spacetime when $\alpha=1$ and $\beta=\frac 12$, with $\lambda$ playing the role of the cosmological contant \cite{CRT}.
\end{itemize}

\noindent
Equation \eqref{GRS} corresponds to an overdetermined system of partial differential equations
of finite type. The study of this system was undertaken in the fundamental paper \cite{NR}. Explicit solutions were determined in \cite{NR} in the two-dimensional case. For the three-dimensional case, the authors restricted in \cite{NR} to the case with $\alpha =0$. Note that as already pointed out in \cite{NR}, for $\alpha=0\neq \beta$, rescaling the vector field $X$ to $-\frac{1}{\beta} X$, equation~\eqref{GRS} reduces to the Ricci soliton equation. We also observe that a trivial solution of \eqref{GRS} is given by $X=0$ and $\beta=\lambda=0$, so we shall always exclude this solution.

The aim of this paper is to determine the three-dimensional homogeneous models of generalized Ricci solitons. A connected, complete and simply connected three-dimensional homogeneous manifold, if not symmetric, is isometric to some Lie group equipped with a left-invariant metric (see \cite{S} for the Riemannian case and \cite{C} for the Lorentzian one). Moreover, with the obvious exceptions of $\mathbb R \times \mathbb S^2$ (Riemannian) and $\mathbb R_1 \times \mathbb S^2$ (Lorentzian), three-dimensional connected simply connected symmetric  spaces are also realized in terms of suitable left-invariant metrics on Lie groups \cite{C2}. For this reason, we shall consider three-dimensional Lie groups, equipped with a left-invariant metric (either Riemannian or Lorentzian). We shall specify our study to solutions of \eqref{GRS} determined by a left-invariant vector field $X$. In this way, \eqref{GRS} will be transformed into a system of algebraic equations, which we can solve, obtaining a complete classification of three-dimensional left-invariant generalized Ricci solitons, and determining several new solutions of \eqref{GRS}. We recall that the study of three-dimensional Ricci solitons already showed some interesting differences arising between the Riemannian case (for which left-invariant solutions do not occur \cite{Cerbo}) and the Lorentzian one, where several left-invariant solutions exist \cite{BCGG}. Also for the broader class of generalized Ricci solitons, interesting differences show up between the Riemannian and the Lorentzian cases. Calculations have been checked by means of {\em Maple 16} $^\copyright$.

\section{$3D$ Riemannian left-invariant generalized Ricci solitons}
\setcounter{equation}{0}

Three-dimensional Riemannian Lie groups were classified in \cite{M}. We shall treat separately the unimodular and non-unimodular cases.

\subsection{Unimodular case}

Let $G$ be a connected three-dimensional Lie group with a left-invariant Riemannian metric. Choose an orientation for the Lie algebra $\mathfrak g$ of $G$, so that the cross product $\times$ is defined on $\mathfrak g$.
The Lie algebra $\mathfrak g$ is unimodular if and only if the endomorphism $L$, defined by $[Z, Y] = L(Z \times Y )$,
is self-adjoint (\cite{M},\cite{HL}). Therefore, $\g$ admits an orthonormal basis $\{e_1,e_2,e_3\}$ of eigenvectors for $L$, so that 
\begin{equation} \label{uniRie}
[e_1, e_2]=C e_3, \qquad  [e_2,e_3]=A e_1, \qquad [e_3,e_1]= B e_2,
\end{equation}
for some real constants $A,B,C$. Explicitly, depending on the sign of $A,B,C$, the Lie group $G$ is isomorphic to one of the cases listed in the following Table~I.

\begin{center}
\begin{tabular}{|c|c|}
\hline  Lie group  & $\vphantom{\frac{{A^A}}{2}}$  $(A,B,C)$ \\
\hline  $\vphantom{\frac{{A^A}}{2}}$ $SU(2)$ & $(+,+,+)$ \\
\hline  $\vphantom{\frac{{A^A}}{2}}$ $\widetilde{SL}(2,\mathbb R)$  & $(+,+,-)$ \\
\hline  $\vphantom{\frac{{A^A}}{2}}$  $\widetilde{E}(2)$ & $(+,+,0)$\\
\hline  $\vphantom{\frac{{A^A}}{2}}$  $E(1,1)$ & $(+,-,0)$\\
\hline  $\vphantom{\frac{{A^A}}{2}}$  $H_3$ & $(+,0,0)$\\
\hline  $\vphantom{\frac{{A^A}}{2}}$  $\mathbb R ^3$ & $(0,0,0)$\\
\hline 
\end{tabular} \nopagebreak \\ \nopagebreak Table I: $3D$ unimodular Riemannian Lie groups $\vphantom{\frac{a}{2}}$
\end{center}

\noindent
In the above Table~I and throughout the paper, $\widetilde{SL}(2,\mathbb R)$ will denote the universal covering of $SL(2,\mathbb R)$, 
$\widetilde{E}(2)$ the universal covering of the group of rigid motions in the Euclidean two-space, $E(1,1)$ the group of rigid motions of
the Minkowski two-space and  $H_3$ the Heisenberg group. 

The description of the Ricci curvature with respect to the basis $\{e_1,e_2,e_3\}$ is well known (see again \cite{M}). We have:
$$Ric=\left(
\begin{array}{ccc}
\frac 12 (A^2-B^2-C^2)+BC & 0 & 0 \\
0 & \frac 12 (B^2-A^2-C^2)+AC & 0 \\
0 & 0 & \frac 12 (C^2-A^2-B^2)+AB
\end{array}
\right).$$
In particular, it is easily seen that the left-invariant metric is Einstein (equivalently, of constant sectional curvature, since we are in dimension three) if and only if either $A=B=C$, $A-B=C=0$, $A=B-C=0$ or $A-C=B=0$.

We now consider an arbitrary vector field $X$, that is, $X=X_i e_i \in \g$, for some real constants $X_1,X_2,X_3$. Then, by \eqref{uniRie}, we get
$$\mathcal L _X g=\left(
\begin{array}{ccc}
0 &  (A-B)X_3  & (C-A)X_2 \\[2pt]
(A-B)X_3 & 0 & (B-C)X_1 \\[2pt]
(C-A)X_2 &  (B-C)X_1 & 0
\end{array}
\right)$$
with respect to the basis $\{e_1,e_2,e_3\}$. Moreover, since  $\{e_1,e_2,e_3\}$ is orthonormal, for any vector  $X=X_i e_i \in \g$  we have $X^\flat \odot X^\flat (e_i,e_j)=X_i X_j$.
Therefore, equation \eqref{GRS} becomes the following system of algebraic equations:

\begin{equation}\label{sys1}
\left\{
\begin{array}{l}
2 \alpha X_1^2-\beta(A^2-B^2-C^2+2BC)=2\lambda, \\[4pt]
2 \alpha X_2^2+\beta(A^2-B^2+C^2-2AC)=2\lambda, \\[4pt]
2 \alpha X_3^2+\beta(A^2+B^2-C^2-2AB)=2\lambda, \\[4pt]
(A-B)X_3+2\alpha X_1X_2 =0, \\[4pt]
(C-A)X_2+2\alpha X_1X_3 =0, \\[4pt]
(B-C)X_1+2\alpha X_2X_3 =0.
\end{array}
\right.
\end{equation}
By standard calculations we obtain the solutions of \eqref{sys1}, proving the following result.

\begin{theorem}\label{GRSuniRie}
Let $\g$ denote a three-dimensional unimodular Riemannian Lie algebra, as described by \eqref{uniRie} with respect to a suitable orthonormal basis $\{e_1,e_2,e_3\}$. Then, up to a renumeration of $e_1,e_2,e_3$, the nontrivial left-invariant generalized Ricci solitons on $\g$ are the following:

\smallskip
(1) $A=B=C$, $\alpha=0$, $\lambda=-\frac{1}{2} \beta A^2$, for all $\beta$ and $X$: when $A=B=C$ (case of constant sectional curvature on $SU(2)$), all vectors in $\g$ are Killing.

\smallskip
(2) $A=B=C$, $\lambda=-\frac{1}{2} \beta A^2$, $X=0$,  for all $\alpha,\beta$: the metric is Einstein.

\smallskip
(3) $A=B-C=0$, $\lambda=0$, $X=0$,  for all $\alpha,\beta$: the metric is flat.

\smallskip
(4) $A\neq B = C$, $\lambda=\frac{1}{2} \beta A(A-2C)$, $X_1=\pm  \sqrt{\frac{\beta A (A-C)}{\alpha}}$, $X_2=X_3=0$, for any $\alpha,\beta$ such that $\alpha \beta A (A-C)>0$. 
\end{theorem}

\smallskip\noindent
Solutions {\em (1)-(3)} are somewhat \lq\lq trivial\rq\rq, as they correspond to cases of metrics of constant sectional curvature and could be expected. On the other hand, by {\em (4)} we have solutions when just two between $A,B,C$ coincide. By Table~I, this yields nontrivial left-invariant generalized Ricci solitons on $SU(2)$, $\widetilde{SL}(2,\mathbb R)$, $\widetilde{E} (2)$, $H_3$. Observe that $\alpha \neq 0$ in case {\em (4)}, while cases {\em (1)-(3)} are Einstein. Thus, no (nontrivial) left-invariant Ricci solitons occur in three-dimensional Riemannian Lie groups, coherently with the results of \cite{Cerbo}.

Finally, we remark that  in solution {\em (4)}, if $\lambda =0$ then $A=2C$ and so, necessarily $\alpha \beta>0$, which  excludes the possibility of solutions of (PS). On the other hand, if (E-W) holds for a three-dimensional manifold, then $\alpha\beta=-1$. Therefore, condition $\alpha \beta A (A-C)>0$ in {\em (4)} yields $A(A-C)<0$. Similarly, in the case of the vacuum near-horizon geometry equation (VN-H) (considered in \cite{NR} both for Riemannian and Lorentzian two-manifolds), we have $\alpha\beta=1$, so that condition $\alpha \beta A (A-C)>0$ in {\em (4)} yields $A(A-C)>0$. 
Taking into account the above Table~I, we then have the following.

\begin{cor}
Three-dimensional Riemannian Lie group $SU(2)$ gives solutions to the special Einstein-Weyl equation {\em (E-W)}. Three-dimensional Riemannian Lie groups $SU(2)$, $\widetilde{SL}(2,\mathbb R)$, $\widetilde{E}(2)$, $H_3$ give solutions to the  vacuum near-horizon geometry equation {\em (VN-H)}. 
\end{cor}

\subsection{Non-unimodular case}

Let now $\mathfrak g$ denote a three-dimensional non-unimodular Riemannian Lie algebra. 
Then, its unimodular kernel $\mathfrak u$ is two-dimensional. Choosing an orthonormal basis $\{ e_1, e_2, e_3 \}$ so that $e_1$ is orthogonal to $\mathfrak u$ and $[e_1, e_2], [e_1 , e_3]$ are mutually orthogonal \cite{M}, the bracket product is described by
\begin{equation} \label{notuniRie}
[e_1, e_2]=A e_2 +B e_3, \quad  [e_1,e_3]=C e_2 + D e_3, \quad [e_2,e_3]=0, \quad A+D\neq 0, \;  AC+BD =0,
\end{equation}
for some real constants $A,B,C,D$.

\noindent
With respect to the basis $\{e_1,e_2,e_3\}$, the Ricci curvature  is described by (see \cite{M})
$$Ric={ \left(
\begin{array}{ccc}
-A^2-\frac 12 B^2-\frac 12 C^2 -D^2 -BC & 0 & 0 \\
0 & -A^2-\frac 12 B^2+\frac 12 C^2 -AD  & 0 \\
0 & 0 & -D^2+\frac 12 B^2-\frac 12 C^2 -AD 
\end{array}
\right).}$$
In particular, the left-invariant metric is of constant sectional curvature if and only if $A-D=B+C=0$.

For  an arbitrary left-invariant vector field $X=X_i e_i \in \g$,  we have
$$\mathcal L _X g=\left(
\begin{array}{ccc}
0 & AX_2 +CX_3   & BX_2 +DX_3  \\[4pt]
AX_2 +CX_3  & -2AX_1 & -(B+C)X_1 \\[4pt]
BX_2 +DX_3  & -(B+C)X_1 & -2D X_1
\end{array}
\right)$$
with respect to the basis $\{e_1,e_2,e_3\}$, and we have again $X^\flat \odot X^\flat (e_i,e_j)=X_i X_j$.
Hence, equation \eqref{GRS} now gives
\begin{equation}\label{sys2}
\left\{
\begin{array}{l}
2 \alpha X_1^2+\beta(2A^2+B^2+C^2 +2D^2 +2BC )=2\lambda, \\[4pt]
-2AX_1+2 \alpha X_2^2+\beta(2A^2+B^2-C^2 +2AD)=2\lambda, \\[4pt]
-2DX_1+2 \alpha X_3^2+\beta(2D^2- B^2+C^2 +2AD )=2\lambda, \\[4pt]
AX_2 +CX_3+2\alpha X_1X_2 =0, \\[4pt]
BX_2 +DX_3+2\alpha X_1X_3 =0, \\[4pt]
-(B+C)X_1+2\alpha X_2X_3 =0.
\end{array}
\right.
\end{equation}
We now solve \eqref{sys2} and list its different solutions, proving the following result. %From now on, we shall not report the solutions corresponding to an Einstein metric ( that is, of constant sectional curvature) with $X=0$.

\begin{theorem}\label{GRSnouniRie}
Let $\g$ denote a three-dimensional non-unimodular Riemannian Lie algebra, as described by \eqref{notuniRie} with respect to a suitable orthonormal basis $\{e_1,e_2,e_3\}$. Then,  the nontrivial left-invariant generalized Ricci solitons on $\g$ are the following:

\smallskip
(1) $A-D=B+C=0$, $\lambda=(2 \alpha^2\beta +\alpha) X_1^2$, $X_1=-\frac{A}{\alpha}$, $X_2=X_3=0$,  for all $\alpha \neq 0$ and $\beta$ 
(constant sectional curvature).

\smallskip
(2) $C=D=0$, $\lambda=\frac{1}{2} \beta (2A^2+B^2)$, $X_1=X_2=0$ and $X_3=\pm  \sqrt{{\frac{\beta(A^2+B^2)}{\alpha}}}$, for all $\alpha$ and $\beta$ satisfying $\alpha\beta >0$.

\smallskip
(3) $A-D=B+C=0$, $\lambda=2\beta A^2$, $X=0$,  for all $\alpha,\beta$: the metric is Einstein.

\smallskip
(4) $B=C=0$,  $\alpha=-\frac{A^2+D^2}{\beta (A+D)^2}\neq 0$, $\lambda=0$, $X_2=X_3=0$, for any $\beta \neq 0$ and $X_1 =\beta (A+D)$.

\smallskip
(5) $A=D$, $B=C=0$, $\lambda=A^2(\frac{1}{\alpha} +2\beta)$, for any $\alpha \neq 0$ and $\beta$, with $X_2=X_3=0$ and $X_1=-\frac{A}{\alpha}$ (constant sectional curvature).  
\end{theorem}

\smallskip\noindent
It is easy to check that the above cases {\em (1)} and {\em (5)} are compatible with (PS), since for $\alpha=1$ and $\beta=-\frac 12$ we have $\lambda=0$ in both  {\em (1)} and {\em (5)}). Moreover, solutions of the form {\em (4)}  compatible with (PS) are a special case of (1),  while  all cases {\em (1), (4), (5)} yield solutions compatible with (E-W). Finally, cases {\em (1)}, {\em (2)} and {\em (5)} are compatible with (VN-H). Thus, we have the following.

\begin{cor}
Three-dimensional non-unimodular Riemannian Lie groups give solutions to the special Einstein-Weyl equation {\em (E-W)}, to the  vacuum near-horizon geometry equation {\em (VN-H)} and (in the case of constant sectional curvature) to the equation {\em (PS)} for a metric projective structure with a skew-symmetric Ricci tensor representative. 
\end{cor}

\section{$3D$ Lorentzian left-invariant unimodular generalized Ricci solitons}
\setcounter{equation}{0}

Let now $\times$ denote the Lorentzian vector product on the Minkowski space $\mathbb R^3 _1$, 
induced by the product of the para-quaternions ($e_1 \times e_2 = −e_3,
e_2 \times e_3 = e_1, e_3 \times e_1 = e_2$, for a pseudo-orthonormal basis ${e_1, e_2, e_3}$, with $e_3$ time-like. The Lie bracket [, ] defines the corresponding Lie algebra $\g$, which is
unimodular if and only if the endomorphism $L$, defined by $[Z, Y] = L(Z \times Y )$,
is self-adjoint \cite{R}. Differently from the Riemannian case, in Lorentzian settings $L$ can assume four different standard forms ({\em Segre types}), giving rise to four classes of three-dimensional unimodular Lorentzian Lie algebras:

\smallskip $\g_1$: $L$ is of Segre type $\{3\}$, that is, its minimal polynomial has a triple root.

\smallskip  $\g_2$: $L$ is of Segre type $\{1 z \bar z \}$, that is, it has two complex conjugate eigenvalues.

\smallskip  $\g_3$: $L$ is of Segre type $\{11,1 \}$, that is, diagonalizable.

\smallskip  $\g_4$: $L$ is of Segre type $\{21\}$, that is, its minimal polynomial has a double root.

\smallskip\noindent
We shall treat these cases separately.

\subsection{Lie algebra $\g_1$}
There exists a pseudo-orthonormal basis $\{e_1,e_2,e_3\}$, with $e_3$ time-like, such that
\begin{eqnarray}
& &\left[e_1,e_2 \right]=A e_1-B e_3, \nonumber \\
\mathfrak{g} _1  : & &\left[ e_1,e_3\right]=-A e_1-B e_2, \label{g1}\\
& & \left[e_2,e_3\right]=B e_1 +A e_2 +A  e_3, \qquad A \neq 0.  \nonumber
\end{eqnarray}
If $B \neq 0$, then $G= \widetilde{SL}(2,\mathbb R)$,
while $G=E(1,1)$ when $B=0$.

The curvature of Lorentzian Lie algebra $\g_1$ was completely determined in \cite{C2}. In particular, with respect to $\{e_i\}$, the Ricci tensor is described by 
$$Ric=\left(
\begin{array}{ccc}
-\frac 12 B^2 &  -AB & AB \\[1pt]
-AB & -2A^2-\frac 12 B^2 & 2 A^2  \\[1pt]
AB & 2A^2 & -2A^2+\frac 12 B^2  
\end{array}
\right),$$
and the left-invariant metric is never Einstein. 

For  a left-invariant vector field $X=X_i e_i \in \g$,  we have
$$\mathcal L _X g=\left(
\begin{array}{ccc}
2A(X_2-X_3) & -AX_1  & AX_1  \\[4pt]
-AX_1 & 2AX_3 & -A (X_2+X_3) \\[4pt]
AX_1 & -A(X_2+X_3) & 2A X_2
\end{array}
\right)$$
with respect to the basis $\{e_1,e_2,e_3\}$, and $X^\flat \odot X^\flat (e_i,e_j)=\varepsilon_i \varepsilon_j X_i X_j$ for all indices $i,j$, where $\varepsilon_1=\varepsilon_2=-\varepsilon_3=1$ corresponds to the causal character of $e_1,e_2,e_3$.
Therefore, equation \eqref{GRS} now becomes
\begin{equation}\label{sys3}
\left\{
\begin{array}{l}
2A(X_2-X_3) + 2 \alpha X_1^2+\beta B^2 =2\lambda, \\[4pt]
2AX_3+2 \alpha X_2^2+\beta(4A^2+B^2)=2\lambda, \\[4pt]
2AX_2+2 \alpha X_3^2+\beta(4A^2- B^2)=-2\lambda, \\[4pt]
-AX_1+2\alpha X_1X_2+2\beta AB=0, \\[4pt]
AX_1-2\alpha X_1X_3-2\beta AB=0, \\[4pt]
-A(X_2+X_3)-2\alpha X_2X_3-4\beta A^2=0.
\end{array}
\right.
\end{equation}
Solving \eqref{sys3}, we obtain the following. %From now on, we shall not report the solutions corresponding to an Einstein metric ( that is, of constant sectional curvature) with $X=0$.

\begin{theorem}\label{GRSg1}
Consider the three-dimensional unimodular Lorentzian Lie algebra $\g_1$, as described by \eqref{g1} with respect to a suitable pseudo-orthonormal basis $\{e_1,e_2,e_3\}$, with $e_3$ time-like. Then, the nontrivial left-invariant generalized Ricci solitons on $\g_1$ are the following:

\smallskip
(1) $\beta=0$, $\lambda=0$, $X_1=0$,  $X_2=X_3= -\frac{A}{\alpha}$, for all $A (\neq 0)$, $B$ and $\alpha \neq 0$.

\smallskip
(2) $\alpha=0$,  $\lambda=\frac{1}{2} \beta B^2$,  $X_1=2\beta B$, $X_2=X_3=-2\beta A$, for all $A (\neq 0)$ and $\beta \neq 0$.

\smallskip
(3) $B=0$, $\lambda=0$, $X_1=0$, $X_2=X_3= \frac{(-1\pm \sqrt{1-8\alpha\beta})A}{2\alpha}$, for all $\alpha,\beta$ satisfying $\alpha\beta \leq \frac 18$.
\end{theorem}

\smallskip\noindent
Solution {\em (2)} in Theorem \ref{GRSg1} corresponds to the existence of Ricci solitons on this class of Lorentzian Lie algebras 
\cite{BCGG}.  Solution {\em (3)}, requiring that $\alpha\beta \leq \frac{1}{8}$ and $\lambda=0$, is incompatible with (VN-H), but compatible with (E-W) and (PS).  These observations yield the following result.

\begin{cor}
Three-dimensional Lorentzian Lie group $E(1,1)$, with Lie algebra described by \eqref{g1}, gives solutions to the special Einstein-Weyl equation {\em (E-W)} and the equation {\em (PS)} for a metric projective structure with a skew-symmetric Ricci tensor representative.  
\end{cor}

\subsection{Lie algebra $\g_2$}
There exists a pseudo-orthonormal basis $\{e_1,e_2,e_3\}$, with $e_3$ time-like, such that
\begin{eqnarray}
& &\left[e_1,e_2 \right]=-C e_2-B e_3, \nonumber \\
\mathfrak{g} _2: & &\left[ e_1,e_3\right]=-B e_2+C e_3, \qquad C \neq 0, \label{g2} \\
& & \left[e_2,e_3\right]=A e_1 .  \nonumber
\end{eqnarray}
In this case, $G=\widetilde{SL}(2,\mathbb R)$ if $A \neq 0$, while $G=E(1,1)$ if $A=0$.
With respect to $\{e_i\}$, the Ricci tensor is given by (see \cite{C2}) 
$$Ric=\left(
\begin{array}{ccc}
-\frac 12 A^2 -2C^2 & 0 & 0 \\
0 & \frac 12 A^2-AB & C(A-2B)  \\
0 & C(A-2B)  & -\frac 12 A^2 +AB  
\end{array}
\right),$$
and the left-invariant metric is never Einstein. 

The Lie derivative $\mathcal L _{X} g$ with respect to a vector $X=X_i e_i \in \g$ is described by 
$$\mathcal L _X g=\left(
\begin{array}{ccc}
0 & -C X_2+(A-B) X_3  & (B-A)X_2 -C X_3  \\[4pt]
-C X_2+(A-B) X_3 & 2CX_1 & 0 \\[4pt]
(B-A)X_2 -C X_3 & 0 & 2C X_1
\end{array}
\right)$$
with respect to the basis $\{e_1,e_2,e_3\}$, and again $X^\flat \odot X^\flat (e_i,e_j)=\varepsilon_i \varepsilon_j X_i X_j$ for all indices $i,j$. Thus, equation \eqref{GRS} becomes
\begin{equation}\label{sys4}
\left\{
\begin{array}{l}
2 \alpha X_1^2+\beta ( A^2 +4C^2)  =2\lambda, \\[4pt]
2CX_1+2 \alpha X_2^2-\beta(A^2-2AB)=2\lambda, \\[4pt]
2CX_1+2 \alpha X_3^2+\beta(A^2 -2AB )=-2\lambda, \\[4pt]
-C X_2+(A-B) X_3+2\alpha X_1X_2=0, \\[4pt]
(B-A)X_2 -C X_3-2\alpha X_1X_3=0, \\[4pt]
-2\alpha X_2X_3-2\beta C(A-2B)=0.
\end{array}
\right.
\end{equation}
We then solve \eqref{sys4}, obtaining the following classification result. %From now on, we shall not report the solutions corresponding to an Einstein metric ( that is, of constant sectional curvature) with $X=0$.

\begin{theorem}\label{GRSg2}
Consider the three-dimensional unimodular Lorentzian Lie algebra $\g_2$, as described by \eqref{g2} with respect to a suitable pseudo-orthonormal basis $\{e_1,e_2,e_3\}$, with $e_3$ time-like. Then, the nontrivial left-invariant generalized Ricci solitons on $\g_2$ are given by 

\smallskip
(1) $A=-2B=\frac{4\alpha X_2 X_3}{3\varepsilon \sqrt{X_3^2 -X_2^2}}$, $C=\varepsilon \alpha\sqrt{X_3^2 -X_2^2}$, 
$\beta=-\frac{3}{8\alpha}$,  $\lambda=\frac{\alpha(3 X_2 ^4-10 X_2^2 X_3 ^2+3 X_3^4)}{ X_2^2- X_3^2\vphantom{a^A}}$, 
$X_1=  -\frac{\varepsilon (X_2^2+X_3^2)}{2\sqrt{X_3^2 -X_2^2}}$, with $\varepsilon=\pm 1$, for all $\alpha \neq 0$.
\end{theorem}

\noindent
We observe that because of condition $\alpha \beta=-\frac{3}{8}$, the above solution is not compatible with any of equations (RS), (E-W), (PS) and (VN-H).

\subsection{Lie algebra $\g_3$}
For a pseudo-orthonormal basis $\{e_1,e_2,e_3\}$, with $e_3$ time-like, of eigenvector of $L$, we have 
\begin{equation}\label{g3}
\mathfrak{g} _3: \quad \left[e_1,e_2 \right]=-C e_3, \quad \left[ e_1,e_3\right]=-B e_2, \quad
 \left[e_2,e_3\right]=A e_1 .
\end{equation}
The following Table~II lists all the Lie groups $G$ which admit a Lie algebra {\bf $\g_3$},
according to the different possibilities for $A$, $B$ and $C$:
\begin{center}
\begin{tabular}{|c|c|}
\hline Lie group &  $(A,B,C)$ $\vphantom{{A^{B^{C}}}}$\\
\hline  $\widetilde{SL}(2,\mathbb R)$ & $(+,+,+)$ $\vphantom{{A^{B^{C^{D}}}}}$\\
\hline  $\widetilde{SL}(2,\mathbb R)$ & $(+,-,-)$  $\vphantom{{A^{B^{C^{D}}}}}$ \\
\hline  $SU(2)$ & $(+,+,-)$ $\vphantom{{A^{B^{C}}}}$\\
\hline  $\widetilde{E}(2)$ & $(+,+,0)$  $\vphantom{{A^{B^{C^{D}}}}}$ \\
\hline  $\widetilde{E}(2)$ & $(+,0,-)$  $\vphantom{{A^{B^{C^{D}}}}}$ \\
\hline  $E(1,1)$ & $(+,-,0)$ $\vphantom{{A^{B^{C}}}}$\\
\hline  $E(1,1)$ & $(+,0,+)$ $\vphantom{{A^{B^{C}}}}$\\
\hline  $H_3$ & $(+,0,0)$  $\vphantom{{A^{B^{C}}}}$\\
\hline  $H_3$ & $(0,0,-)$ $\vphantom{{A^{B^{C}}}}$\\
\hline  $\mathbb R \oplus \mathbb R \oplus \mathbb R$ & $(0,0,0)$ $\vphantom{{A^{B^{C}}}}$\\
\hline
\end{tabular} \nopagebreak \\ \nopagebreak {\em Table~II: $3D$ Lorentzian Lie groups with Lie algebra $\g_3$} $\vphantom{\frac{a}{2}}$
\end{center}

\noindent
Following \cite{C2} (or by direct calculation), the Ricci curvature of Lorentzian Lie algebra $\g_3$, with respect to $\{e_i\}$, is described by 
$$Ric=\left(
\begin{array}{ccc}
\frac 12 (B^2-A^2+C^2)-BC & 0 & 0 \\
0 & \frac 12 (A^2- B^2+C^2)-AC & 0 \\
0& 0 & \frac 12 (C^2-A^2-B^2)+AB  
\end{array}
\right),$$
and the left-invariant metric is Einstein (equivalently, of constant sectional curvature) if and only if either $A=B=C$, $A-B=C=0$, $A-C=B=0$ or $A=B-C=0$. In the last three cases, the metric is flat.  

For  a left-invariant vector field $X=X_i e_i \in \g$,  we have
$$\mathcal L _X g=\left(
\begin{array}{ccc}
0 &  (A-B) X_3  & (C-A) X_2 \\[2pt]
(A-B) X_3 & 0 & (B-C) X_1 \\[2pt]
 (C-A) X_2 & (B-C) X_1 & 0
\end{array}
\right)$$
with respect to the basis $\{e_1,e_2,e_3\}$. Finally, $X^\flat \odot X^\flat (e_i,e_j)=\varepsilon_i \varepsilon_j X_i X_j$ for all indices $i,j$. Hence, equation \eqref{GRS} now yields
\begin{equation}\label{sys5}
\left\{
\begin{array}{l}
 2 \alpha X_1^2+\beta (A^2- B^2-C^2+2BC) =2\lambda, \\[4pt]
2 \alpha X_2^2-\beta(A^2-B^2+C^2-2AC)=2\lambda, \\[4pt]
2 \alpha X_3^2+\beta(A^2+B^2-C^2-2AB)=-2\lambda, \\[4pt]
(A-B)X_3+2\alpha X_1X_2 =0, \\[4pt]
(C-A)X_2-2\alpha X_1X_3 =0, \\[4pt]
(B-C)X_1-2\alpha X_2X_3 =0.
\end{array}
\right.
\end{equation}
Although system \eqref{sys5} is rather similar to \eqref{sys1}, the different signs, due to the different signature of the metric, are responsible for the existence of many more solutions for \eqref{sys5}. Solving \eqref{sys5}, we prove the following. %From now on, we shall not report the solutions corresponding to an Einstein metric ( that is, of constant sectional curvature) with $X=0$.

\begin{theorem}\label{GRSg3}
Consider the three-dimensional unimodular Lorentzian Lie algebra $\g_3$, as described by \eqref{g3} with respect to a suitable pseudo-orthonormal basis $\{e_1,e_2,e_3\}$, with $e_3$ time-like. Then, the nontrivial left-invariant generalized Ricci solitons on $\g_3$ are the following:

\smallskip
(1) $A=B=C$, $\alpha=0$, $\lambda=\frac 12 \beta A^2$, for all $X$: when $A=B=C$, all left-invariant vector fields are Killing.

\smallskip
(2) $A=B-C=0$,  $\alpha=\lambda=0$,  $X_2=X_3=0$: $X=X_1 e_1$ is Killing for the flat metric obtained when $A=B-C=0$ (corresponding solutions occur for $A-B=C=0$ and $A-C=B=0$). 

\smallskip
(3) $A=B-C=0$, $\lambda=0$, $X=0$,  for all $\alpha,\beta$: the metric is flat (corresponding solutions occur for $A-B=C=0$ and $A-C=B=0$).

\smallskip
(4) $A=B=C$, $\lambda=\frac{1}{2} \beta A^2$, $X=0$,  for all $\alpha,\beta$: the metric is Einstein.

\smallskip
(5) $B=C\neq A$, $\lambda=\frac{1}{2} \beta A (2C-A)$, $X_1=\pm \sqrt{\frac{\beta A(C-A)}{\alpha}}$, $X_2=X_3=0$,  for all $\alpha,\beta$, with $\alpha \beta A(C-A) >0$.

\smallskip
(6) $A=B\neq C$, $\lambda=\frac{1}{2} \beta C (2B-C)$, $X_1=X_2=0$, $X_3=\pm \sqrt{\frac{\beta C(C-B)}{\alpha}}$, for all $\alpha,\beta$, with $\alpha \beta C(C-B) >0$.

\smallskip
(7) $A+B=C=0$, $\beta=-\frac{3}{8\alpha}$, $\lambda=\frac{A^2}{2\alpha}$,  $X_1=-X_2= \frac{\varepsilon A}{\sqrt{2} \alpha}$, $X_3= \frac{\varepsilon A}{2\alpha}$, $\varepsilon=\pm 1$,  for all $\alpha \neq 0$.

\smallskip
(8) $A=-\frac{\varepsilon(2 X_1^2+X_2^2)}{4\beta\sqrt{X_1^2+X_2^2}}$, $B=\frac{\varepsilon(X_1^2+2X_2^2)}{4\beta\sqrt{X_1^2+X_2^2}}$, 
$C= \frac{\varepsilon(X_1^2-X_2^2)}{4\beta \sqrt{X_1^2+X_2^2}}$, $\alpha=-\frac{3}{8\beta}$, $\lambda=-\frac{X_1^4+X_1^2X_2^2+X_2^4}{4\beta (X_1^2+X_2^2)}$, $X_3=-\frac{\varepsilon X_1 X_2}{\sqrt{X_1^2+X_2^2}}$, for all $\beta \neq 0$.
\end{theorem}

\smallskip\noindent
%With regard to solution {\em (8)}, we observe that equation \eqref{X1bar} admits real solutions for all values of $A,B$. 
Since $\alpha \beta=-\frac{3}{8}$, solutions {\em (7)} and {\em (8)} are not compatible with any of equations (RS), (E-W), (PS) and (VN-H). We now focus our attention to solutions {\em (5)} and {\em (6)}. Note that these solutions are very similar to one another. However, we listed both of them, since for {\em (5)} the vector satisfying \eqref{sys5} is space-like, while for {\em (6)} is time-like. %Moreover, we observe that conditions $A=B$ and $B=C$ are also related to the classification of three-dimensional naturally reductive Lorentzian Lie groups \cite{CM}.

Solutions (5) and (6) are compatible with (E-W), (PS) and (VN-H). More precisely, for any choice of $A$ and $C\neq A$:
\begin{itemize}
\item If $A(A-C)>0$, then by {\em (5)} we have that the same left-invariant Lorentzian metric on $\g_3$ is solution to both (E-W) and (PS).
\item If $A(A-C)<0$, then {\em (5)} describes a left-invariant Lorentzian metric on $\g_3$, which is a solution to (VN-H).
\end{itemize}
Similar observations hold for solution {\em (6)}, discussing the cases when $\alpha \beta C(C-B) >0$. Hence, taking into account Table~II, we have the following result.

\begin{cor}
Three-dimensional Lorentzian Lie groups $\widetilde{SL}(2,\mathbb R)$ and $H_3$, with Lie algebra described by \eqref{g3}, give solutions to the special Einstein-Weyl equation {\em (E-W)},  the equation {\em (PS)} for a metric projective structure with a skew-symmetric Ricci tensor representative and to the  vacuum near-horizon geometry equation {\em (VN-H)}.  Three-dimensional Lorentzian Lie group $SU(2)$, with Lie algebra described by \eqref{g3}, gives solutions to {\em (VN-H)}.
\end{cor}

\subsection{Lie algebra $\g_4$}
There exists a pseudo-orthonormal basis $\{e_1,e_2,e_3\}$, with $e_3$ time-like, such that
\begin{eqnarray}
& &\left[e_1,e_2 \right]=- e_2 + (2 \eta - B) e_3, \qquad \eta = \pm 1, \nonumber \\
\mathfrak{g} _4: & &\left[ e_1,e_3\right]=-B e_2 + e_3,  \label{g4} \\
& & \left[e_2,e_3\right]=A e_1 .  \nonumber
\end{eqnarray}
Table~III below describes all Lie groups $G$ admitting a Lie algebra {\bf $g_4$}:
\begin{gather*}
\begin{array}{c}
\begin{tabular}{|c|c|c|}
\hline $\quad$ Lie group $\quad$ & $\quad$ $\eta A$ $\quad$ & $\quad$ $B$ $\quad$ $\vphantom{{A^{B^{C}}}}$ \\
\hline  $\widetilde{SL}(2,\mathbb R)$ & $\neq 0$ & $\neq \eta$ $\vphantom{{A^{B^{C^{D}}}}}$ \\
\hline  $E(1,1)$ & $0$ & $\neq \eta$ $\vphantom{{A^{B^{C}}}}$ \\
\hline  $E(1,1)$ & $<0$ & $\eta$  $\vphantom{{A^{B^{C}}}}$ \\
\hline  $\widetilde{E}(2)$ & $>0$ & $\eta$ $\vphantom{{A^{B^{C^{D}}}}}$ \\
\hline  $H_3$ & $0$ & $\eta$ $\vphantom{{A^{B^{C}}}}$ \\
\hline
\end{tabular}
\end{array}
 \\
 \text{{\em Table~III: $3D$ Lorentzian Lie groups with Lie algebra $\g_4$}}
 \vphantom{\frac{a}{2}}
 \end{gather*}

With respect to $\{e_i\}$, the Ricci tensor is described by (see \cite{C2})
$$Ric=\left(
\begin{array}{ccc}
-\frac 12 A^2 & 0 & 0 \\
0 & \frac 12 A^2+2\eta (A-B) -AB+2 & A+2(\eta-B)  \\
0 & A+2(\eta-B) & -\frac 12 A^2+AB +2-2\eta B  
\end{array}
\right),$$
and the left-invariant metric is Einstein if and only if $A=B-\eta=0$. 

For a left-invariant vector field $X=X_i e_i \in \g$,  we have
$$\mathcal L _X g=\left(
\begin{array}{ccc}
0 & -X_2+(A-B)X_3  & (B-A-2\eta)X_2 -X_3  \\[4pt]
-X_2+(A-B)X_3 & 2X_1 & 2\eta X_1 \\[4pt]
(B-A-2\eta)X_2-X_3 & 2\eta X_1 & 2X_1
\end{array}
\right)$$
with respect to the basis $\{e_1,e_2,e_3\}$, and $X^\flat \odot X^\flat (e_i,e_j)=\varepsilon_i \varepsilon_j X_i X_j$.
Therefore, equation \eqref{GRS} now becomes
\begin{equation}\label{sys6}
\left\{
\begin{array}{l}
 2 \alpha X_1^2+\beta A^2 =2\lambda, \\[4pt]
2 X_1+2 \alpha X_2^2-\beta(A^2+4\eta (A-B) -2AB+4)=2\lambda, \\[4pt]
2 X_1+2 \alpha X_3^2+\beta(A^2-2AB-4+4\eta B )=-2\lambda, \\[4pt]
-X_2+(A-B)X_3+2\alpha X_1X_2=0, \\[4pt]
(B-A-2\eta)X_2 -X_3-2\alpha X_1X_3=0, \\[4pt]
2\eta X_1-2\alpha X_2X_3-2\beta ( A+2(\eta-B))=0.
\end{array}
\right.
\end{equation}
We then solve \eqref{sys6} and prove the following. %From now on, we shall not report the solutions corresponding to an Einstein metric ( that is, of constant sectional curvature) with $X=0$.

\begin{theorem}\label{GRSg4}
Consider the three-dimensional unimodular Lorentzian Lie algebra $\g_4$, as described by \eqref{g4} with respect to a suitable pseudo-orthonormal basis $\{e_1,e_2,e_3\}$, with $e_3$ time-like. Then, the nontrivial left-invariant generalized Ricci solitons on $\g_4$ are the following:

\smallskip
(1) $A=B-\eta$, $\lambda=\frac{1}{2}\beta A^2$, $X_1=0$,  $X_2=-\eta X_3$,  
$X_3=  \pm\sqrt{-\frac{\eta \beta A}{\alpha}}$, for all $\alpha, \beta$ satisfying $\eta A \alpha\beta <0$.

\smallskip
(2) $A=B-\eta$, $\alpha=0$, $\lambda=\frac{1}{2}\beta A^2$, $X_1=-\eta \beta A$,  $X_2=-\eta X_3$, for any value of $\beta$ and 
$X_3$.

\smallskip
(3) $B=\frac 12 A+\eta$, $\beta=-\frac{1}{8\alpha}$, $\lambda=0$, $X_1=\frac{\eta A}{4\alpha}$,  
$X_2= \pm \sqrt{-\frac{\eta A}{4\alpha^2}}$, $X_3=-\eta X_2$,  for any $\alpha \neq 0$ and $\eta A <0$. 

\medskip
(4) $\beta= -\frac{A(A-B+\eta)}{\alpha (A-2B+2\eta)^2}$, $\lambda= - \frac 12 \beta A (A-2B+2\eta)$, $X_1=\beta(A-2B+2\eta)$,  $X_2=X_3=0$,  for any $\alpha \neq 0$ and $A-2B+2\eta \neq 0$. 

\medskip
(5) $\beta= -\frac{A-B+\eta}{4\alpha A}$, $\lambda=  \frac{(A-B+\eta) (A-2B+2\eta)}{8\alpha}$, 
$X_1= \frac{A-B+\eta}{2\eta \alpha }$,  $X_2=-\eta X_3$ and  
$$X_3 = \pm \frac{1}{2\alpha}\sqrt{\frac{5\eta AB -3\eta A^2-5A-2\eta+4B-2\eta B^2}{A}},$$
for any $\alpha \neq 0$ and $A \neq 0$. 

\end{theorem}

\smallskip\noindent
In the above Theorem~\ref{GRSg4}, solution {\em (2)} corresponds to the existence of Ricci solitons for Lorentzian Lie algebras of the form $\g_4$ \cite{BCGG}.  Solution {\em (1)} requires that the sign of  $\alpha\beta$ is opposite to the one of $\eta A$. Consequently, taking into account the above Table~III,  it yields solutions to (E-W) for Lorentzian Lie groups $\widetilde{SL}(2,\mathbb R)$ and $\widetilde{E}(2)$, and solutions to (VN-H) for $\widetilde{SL}(2,\mathbb R)$ and ${E}(1,1)$. 

In solution {\em (4)}, we have  $\alpha\beta= -\frac{A(A-B+\eta)}{(A-2B+2\eta)^2}$. For (E-W), as $\alpha\beta=-1$, this equation becomes $4B^2+(5A+8\eta)B+3\eta A=0$, which admits real solutions for any value of $A$. On the other hand, for (VN-H) we have $\alpha\beta=\frac 12$. Hence, we find $3A^2-6(B-\eta)A+4(B-\eta)^2=0$, whih has not real solutions. Finally, it is easily seen that solution {\em (4)} is not compatible with (PS).

By similar arguments, from solution {\em (5)} we find again  $\widetilde{SL}(2,\mathbb R)$ as solution to (E-W) and (VN-H). Thus, we proved the following.

\begin{cor}
Three-dimensional Lorentzian Lie group  $\widetilde{SL}(2,\mathbb R)$, with Lie algebra described by \eqref{g4}, gives solutions to the special Einstein-Weyl equation {\em (E-W)} and the vacuum near-horizon equation {\em (VN-H)}. Moreover, Lorentzian Lie groups $\widetilde{E}(2)$ and $E(1,1)$, with Lie algebra described by \eqref{g4}, give solutions to {\em (E-W)} and {\em (VN-H)} respectively. 
 
\end{cor}

\section{$3D$ Lorentzian left-invariant non-unimodular generalized Ricci solitons}
\setcounter{equation}{0}

Let now $\mathfrak g$ denote a three-dimensional non-unimodular Lorentzian Lie algebra. Differently from the Riemannian case, we must now consider three distinct cases, depending on whether the two-dimensional unimodular kernel $\mathfrak u$ is either space-like, time-like or degenerate \cite{CoPa}. These three cases were listed in \cite{C} as Lorentzian Lie algebras $\g_5$, $\g_6$ and $\g_7$.

\subsection{Lie algebra $\g_5$}
There exists a pseudo-orthonormal basis $\{ e_1, e_2, e_3 \}$, with $e_3$ time-like, such that
\begin{equation} \label{g5}
[e_1, e_2]=0, \quad  [e_1,e_3]=A e_1 + B e_2, \quad [e_2,e_3]=C e_1+D e_2, \quad A+D\neq 0, \;  AC+BD =0,
\end{equation}
for some real constants $A,B,C,D$.

\noindent
With respect to the basis $\{e_1,e_2,e_3\}$, the Ricci curvature  is described by (see \cite{C2})
$$Ric={\left(
\begin{array}{ccc}
A^2+\frac 12 B^2-\frac 12 C^2+AD & 0 & 0 \\
0 & AD-\frac 12 B^2+\frac 12 C^2 +D^2  & 0 \\
0 & 0 & -A^2-\frac 12 B^2-\frac 12 C^2 -D^2-BC 
\end{array}
\right)}.$$
In particular,the left-invariant metric is of constant sectional curvature if and only if $A-D=B+C=0$.

For  an arbitrary left-invariant vector field $X=X_i e_i \in \g$,  we have
$$\mathcal L _X g=\left(
\begin{array}{ccc}
2A X_3 & (B+C) X_3   & -A X_1 -C X_2  \\[4pt]
(B+C) X_3  & 2D X_3 & -B X_1 -D X_2 \\[4pt]
-A X_1 -C X_2  & -B X_1 -D X_2 & 0
\end{array}
\right).$$
With respect to the basis $\{e_1,e_2,e_3\}$, we have again $X^\flat \odot X^\flat (e_i,e_j)=\varepsilon_i \varepsilon_jX_i X_j$.
Hence, equation \eqref{GRS} now gives
\begin{equation}\label{sys7}
\left\{
\begin{array}{l}
2A X_3 +2 \alpha X_1^2-\beta(2A^2+ B^2- C^2+2AD)=2\lambda, \\[4pt]
2D X_3+2 \alpha X_2^2-\beta(2AD- B^2+ C^2 +2D^2)=2\lambda, \\[4pt]
2 \alpha X_3^2+\beta(2 A^2+ B^2+ C^2+2D^2+2 BC )=-2\lambda, \\[4pt]
(B+C)X_3+2\alpha X_1X_2 =0, \\[4pt]
-A X_1 -CX_2 -2\alpha X_1X_3 =0, \\[4pt]
-B X_1 -D X_2-2\alpha X_2X_3 =0.
\end{array}
\right.
\end{equation}
We then solve \eqref{sys7} and obtain the following.

\begin{theorem}\label{GRSg5}
Let $\g$ denote a three-dimensional non-unimodular Lorentzian Lie algebra $\g_5$, as described by \eqref{g5} with respect to a suitable pseudo-orthonormal basis $\{e_1,e_2,e_3\}$, with $e_3$ time-like. Then, the nontrivial left-invariant generalized Ricci solitons on $\g_5$ are the following:

\smallskip
(1) $A-D=B+C=0$, $\lambda=-(\frac{1}{\alpha}+2\beta)A^2$, $X_1=X_2=0$, $X_3=-\frac{A}{\alpha}$,   for all $\alpha \neq 0$ and $\beta$ 
(constant sectional curvature).

\smallskip
(2) $C=D=0$, $\beta=-\frac{1}{4\alpha}$,  $\lambda=\frac{B^2}{8\alpha}$, $X_1=\varepsilon X_3$, $X_2=\frac{\varepsilon B}{2\alpha}$, 
$X_3 =-\frac{A}{2\alpha}$, $\varepsilon = \pm 1$, for all $\alpha \neq 0$.

\smallskip
(3) $A-D=B+C=0$, $\lambda=-2\beta A^2$, $X=0$,  for all $\alpha,\beta$: the metric is Einstein.

\smallskip
(4) $C=D=0$, $\lambda=-\frac{1}{2}\beta(2A^2+B^2)$, $X_1=X_3=0$, $X_2 = \pm  \sqrt{-\frac{\beta(A^2+B^2)}{\alpha}}$, for any $\alpha$ and $\beta$ satisfying $\alpha \beta <0$.

\smallskip
(5) $B=C=D=0$, $\beta=-\frac{1}{\alpha}$, $\lambda=0$, $X_1=X_2=0$, $X_3 = - \frac{A}{\alpha}$, for any $\alpha \neq 0$.

%\smallskip
%(6) $A-D=B=C=0$,  $\lambda=-(\frac{1}{\alpha}+2\beta)D^2$, $X_1=X_2=0$, $X_3 = - \frac{D}{\alpha}$, for any $\alpha \neq 0$ and $\beta$.

\medskip
(6) $B=C=0$,  $\lambda= \frac{A^2(32\alpha ^3 \beta^3 +28 \alpha^2 \beta ^2 +9\alpha \beta+1)}{4\alpha (2\alpha\beta+1)^2 }$, $X_1= \pm \frac{A\sqrt{8\alpha^2 \beta^2+5\alpha\beta+1}}{2\alpha (2\alpha\beta+1)}$, $X_2=0$, $X_3 = - \frac{A}{2\alpha}$, for any $\alpha \neq 0$ and $\beta$, with $2\alpha\beta+1 \neq 0$ (note that $a^2 \beta^2+5\alpha\beta+1>0$).

\medskip
(7) $B=C=0$, $\beta=- \frac{2A-D}{4\alpha(A-D)}$,  $\lambda=- \frac{A(2A^2-AD+D^2)}{4\alpha(A-D)}$, $X_1=0$,  
$X_2=\pm \frac{1}{\alpha}\sqrt{A^2-\frac 12 AD+\frac 12 D^2}$, $X_3 = - \frac{D}{2\alpha}$, for any $\alpha \neq 0$ and $A\neq D$ (note that $2A^2-AD+D^2>0$).

%\smallskip
%(9) $A-D=B-C=0$, $\beta=-\frac{3}{8\alpha}$, $X_2=-kX_1$, $X_3 =k\alpha X_1^2$, where $k:=\pm\sqrt{\frac{1}{4\alpha^2 X_1^2-C^2}}$ and
%$X_1$ is a real solution (when it exists)  of
%
%$$2\alpha^3 X_1^4-4\alpha^2 \lambda X_1^2+\lambda C^2=0.$$
%
\end{theorem}

\smallskip\noindent
It is easy to check that several of the above solutions are compatible with (E-W),(PS) and (VN-H). In particular, solution {\em (1)} yields solutions to all these equations. Hence, we have the following.

\begin{cor}
Three-dimensional non-unimodular Lorentzian Lie groups with Lie algebra $\g_5$ give solutions to the special Einstein-Weyl equation {\em (E-W)}, to the equation {\em (PS)} for a metric projective structure with a skew-symmetric Ricci tensor representative and to the vacuum near-horizon geometry equation {\em (VN-H)}. 
\end{cor}

\subsection{Lie algebra $\g_6$}
There exists a pseudo-orthonormal basis $\{ e_1, e_2, e_3 \}$, with $e_3$ time-like, such that
\begin{equation} \label{g6}
[e_1, e_2]=A e_2 + B e_3, \quad  [e_1,e_3]=C e_2+D e_3, \quad [e_2,e_3]=0, \quad A+D\neq 0, \;  AC-BD =0,
\end{equation}
for some real constants $A,B,C,D$.

\noindent
Following  \cite{C2}, with respect to the basis $\{e_1,e_2,e_3\}$ we have
$$Ric={\left(
\begin{array}{ccc}
\frac 12 B^2-A^2+\frac 12 C^2-D^2-BC & 0 & 0 \\
0 & \frac 12 B^2-A^2-\frac 12 C^2 -AD  & 0 \\
0 & 0 & AD+\frac 12 B^2-\frac 12 C^2 +D^2 
\end{array}
\right)}.$$
In particular, the left-invariant metric is of constant sectional curvature if and only if either $A-D=B-C=0$, $A+B=C+D=0$ or $A-B=C-D=0$.

For  an arbitrary left-invariant vector field $X=X_i e_i \in \g$,  we have
$$\mathcal L _X g=\left(
\begin{array}{ccc}
0 & A X_2 +C X_3   & -B X_2 -D X_3  \\[4pt]
A X_2 +C X_3  & -2A X_1 & (B-C) X_1 \\[4pt]
-B X_2 -D X_3  & (B-C)X_1 & 2D X_1
\end{array}
\right).$$
Equation \eqref{GRS} now becomes
\begin{equation}\label{sys8}
\left\{
\begin{array}{l}
2 \alpha X_1^2+\beta(2A^2- B^2- C^2+2D^2+2BC)=2\lambda, \\[4pt]
-2A X_1+2 \alpha X_2^2+\beta(2A^2- B^2+ C^2+2AD)=2\lambda, \\[4pt]
2D X_1+2 \alpha X_3^2-\beta( 2AD+ B^2- C^2 +2D^2)=-2\lambda, \\[4pt]
A X_2 +C X_3  +2\alpha X_1X_2 =0, \\[4pt]
-B X_2 -D X_3 -2\alpha X_1X_3 =0, \\[4pt]
(B-C) X_1-2\alpha X_2X_3 =0.
\end{array}
\right.
\end{equation}
We solve \eqref{sys8} and prove the following.

\begin{theorem}\label{GRSg6}
Let $\g$ denote a three-dimensional non-unimodular Lorentzian Lie algebra $\g_6$, as described by \eqref{g6} with respect to a suitable pseudo-orthonormal basis $\{e_1,e_2,e_3\}$, with $e_3$ time-like. Then, the nontrivial left-invariant generalized Ricci solitons on $\g_6$ are the following:

\smallskip
(1) $A-D=B-C=0$, $\lambda=(\frac{1}{\alpha}+2\beta)A^2$, $X_1=-\frac{A}{\alpha}$, $X_2=X_3=0$,  for all $\alpha \neq 0$ and $\beta$ 
(constant sectional curvature).

\smallskip
(2) $A-D=B-C=0$, $\lambda=2\beta A^2$, $X=0$,  for all $\alpha$ and $\beta$ 
(Einstein).

\smallskip
(3) $B=\pm A, C=\pm D$, $\lambda=\frac 12 \beta (A+D)^2$, $X=0$,  for all $\alpha,\beta$ (Einstein).

\smallskip
(4) $C=D=0$, $\lambda=\frac{1}{2}\beta(2A^2-B^2)$, $X_1=X_2=0$, $X_3 = \pm  \sqrt{\frac{\beta(B^2-A^2)}{\alpha}}$, for any $\alpha$ and $\beta$ satisfying $\alpha \beta(B^2-A^2) >0$.

%\smallskip
%(5) $B=C=D=0$, $\beta=-\frac{1}{\alpha}$, $\lambda=0$,  $X_1 = - \frac{A}{\alpha}$, $X_2=X_3=0$, for any $\alpha \neq 0$.

\smallskip
(5) $B=C=0$, $\beta=-\frac{A^2+D^2}{\alpha(A+D)^2}$, $\lambda=0$,  $X_1 =-\frac{A^2+D^2}{\alpha(A+D)}$, $X_2=X_3=0$, for any $\alpha \neq 0$.

\medskip
(6) $B=C=0$, $\beta=- \frac{2A-D}{4\alpha(A-D)}$,  $\lambda=- \frac{A(2A^2-AD+D^2)}{4\alpha(A-D)}$, $X_1=- \frac{D}{2\alpha}$,  $X_2=0$, 
$X_3=\pm \frac{1}{2\alpha}\sqrt{2A^2- AD+ D^2}$, $X_3 = 0$, for any $\alpha \neq 0$ and $A\neq D$ (note that $2A^2-AD+D^2>0$).

%\medskip
%(7) $A=B=C=0$, $\beta =-\frac{1}{4\alpha}$,  $\lambda=0$, $X_1=-\frac{D}{2\alpha}$, $X_2=0$, $X_3 = \pm X_1$, for any $\alpha \neq 0$ and $\beta$.

\medskip
(7) $A=B=0$, $\lambda=\frac 12 \beta (2D^2-C^2)$, $X_1=0$, $X_2 = \pm  \sqrt{\frac{\beta(D^2-C^2)}{\alpha}}$, $X_3=0$, for all $\alpha,\beta$ with $\alpha \beta(D^2-C^2)>0$.
\end{theorem}

\smallskip\noindent
By the same argument used in the previous case, it is easy to check that several of the above solutions are compatible with (E-W),(PS) and (VN-H). Thus, we have the following.

\begin{cor}
Three-dimensional non-unimodular Lorentzian Lie groups with Lie algebra $\g_6$ give solutions to the special Einstein-Weyl equation {\em (E-W)}, to the equation {\em (PS)} for a metric projective structure with a skew-symmetric Ricci tensor representative and to the vacuum near-horizon geometry equation {\em (VN-H)}. 
\end{cor}

\subsection{Lie algebra $\g_7$}
There exists a pseudo-orthonormal basis $\{ e_1, e_2, e_3 \}$, with $e_3$ time-like, such that
\begin{equation} \label{g7}
[e_1, e_2]=-[e_1,e_3]=-A e_1 - B e_2-B e_3, \quad  [e_2,e_3]=C e_1 +D e_2+D e_3, \quad A+D\neq 0, \;  AC =0,
\end{equation}
for some real constants $A,B,C,D$. The Ricci curvature is then described by (see \cite{C2})
$$Ric=\left(
\begin{array}{ccc}
-\frac 12 C^2 & 0 & 0 \\[2pt]
0 & AD-A^2+\frac 12 C^2 -BC  & A^2-AD+BC \\[2pt]
0 & A^2-AD+BC & AD-A^2-\frac 12 C^2-BC 
\end{array}
\right)$$
and the left-invariant metric is of constant sectional curvature (flat) if and only if either $A=C=0$ or $A-D=C=0$.
For a vector $X=X_i e_i \in \g$,  we have
$$\mathcal L _X g=\left(
\begin{array}{ccc}
-2A (X_2-X_3) &AX_1-BX_2+(B+C)X_3    & -AX_1+(B-C)X_2-B X_3 \\[4pt]
AX_1-BX_2+(B+C)X_3  & 2B X_1+2D X_3 & -2B X_1-D X_2-DX_3 \\[4pt]
-AX_1+(B-C)X_2-B X_3  & -2B X_1-D X_2-DX_3 & 2B X_1+2D X_2
\end{array}
\right).$$
Thus, equation \eqref{GRS} becomes
\begin{equation}\label{sys9}
\left\{
\begin{array}{l}
-2A (X_2-X_3) +2 \alpha X_1^2+\beta C^2=2\lambda, \\[4pt]
2B X_1+2D X_3+2 \alpha X_2^2+\beta(2A^2-2AD- C^2 +2BC)=2\lambda, \\[4pt]
2B X_1+2D X_2+2 \alpha X_3^2+\beta(2A^2-2AD+ C^2+2BC)=-2\lambda, \\[4pt]
AX_1-BX_2+(B+C)X_3 +2\alpha X_1X_2 =0, \\[4pt]
-AX_1+(B-C)X_2-B X_3  -2\alpha X_1X_3 =0, \\[4pt]
-2B X_1-D X_2-DX_3-2\alpha X_2X_3-2\beta (A^2-AD+BC) =0.
\end{array}
\right.
\end{equation}
Solving \eqref{sys9}, we prove the following.

\begin{theorem}\label{GRSg7}
Let $\g$ denote a three-dimensional non-unimodular Lorentzian Lie algebra $\g_7$, as described by \eqref{g7} with respect to a suitable pseudo-orthonormal basis $\{e_1,e_2,e_3\}$, with $e_3$ time-like. Then, the nontrivial left-invariant generalized Ricci solitons on $\g_7$ are the following:

\smallskip
(1) $A=C=0$, $\alpha=0$,  $\lambda=0$, $X_2=X_3=-\frac{B}{D}X_1$, for all $\beta$ and $X_1$ (flat).

\smallskip
(2) $A=\frac 12 D $, $C=0$, $\alpha=0$,  $ X_1=-\frac{4B \lambda}{D^2}$, $ X_2=\frac{16 \lambda B^2 -4\lambda D^2+\beta D^4 }{4D^3}$, $ X_3=\frac{16 \lambda B^2 +4\lambda D^2+\beta D^4 }{4D^3}$,  for all $\beta$ and $\lambda$.

\smallskip
(3) $C=0$, $\lambda=0$, $X_1=0$, $X_2=X_3$ solution of
$$\alpha x^2+D x +\beta A(A-D)=0,$$
for all $A,B,D$ and $\alpha,\beta$ such that $D^2-4\alpha \beta A(A-D)>0$.

\smallskip
(4) either $A=C=0$ or $A-D=C=0$, $\lambda=0$, $X=0$: flat metric.

%\smallskip
%(5) $A=\frac 12 D$, $B=C=0$, $\alpha=0$, $X_1 =0$, $X_2=-\frac{\lambda}{D}+\frac 14 \beta D$, $X_3=\frac{\lambda}{D}+\frac 14 \beta D$, for all $\beta$ and $\lambda$.

%\smallskip
%(6) $B=C=0$, $\beta=-\frac{A^2+D^2}{\alpha(A+D)^2}$, $\lambda=0$,  $X_1 =\beta(A+D)$, $X_2=X_3=0$, for any $\alpha \neq 0$.

\medskip
(5) $A=\frac 12 D$, $B=C=0$, $\alpha=0$ $\lambda=\beta A^2$,   $X_1=X_2=0$, $X_3=\beta A$, for any $\beta \neq 0$.

\medskip
(6) $A=0$, $\beta=-\frac{1}{4\alpha}$ $\lambda =\frac{C^2}{8\alpha}$, $X_1=-\frac{C}{2\alpha}$, $X_2=X_3$ solution of
$$4\alpha^2 x^2+4\alpha D x-3B C=0,$$
for any $\alpha \neq 0$, whenever $D^2+3BC>0$.

\medskip
(7) $A=B=0$, $\lambda= -\frac 12 \beta C^2$, $X_1 =  \pm \sqrt{\frac{-\beta C^2}{\alpha}}$, $X_2=X_3=0$, for any $\alpha,\beta$ with $\alpha \beta C^2<0$.

\medskip
(8) $A=0$,  $\beta=-\frac{1}{\alpha}$, $\lambda=-\frac 12 \beta C^2$, $X_1=-\beta C $, $X_2=X_3=0$, for any $\alpha \neq 0$.
\end{theorem}

\smallskip\noindent
The above solutions with $\alpha=0$ correspond to the existence of left-invariant Ricci solitons on $\g_7$ \cite{BCGG}. Moreover, several of the above solutions are compatible with (E-W), (PS) and (VN-H). In particular, whatever the value of $\alpha$ and $\beta$, if $D$ is sufficiently big then $D^2-4\alpha \beta A(A-D)>0$. Henceforth, {\em (3)} yields solutions to all (E-W), (PS) and (VN-H). %Solutions {\em (10)}  are compatible with (E-W) and (PS), respectively.   
Thus, we have the following.

\begin{cor}
Three-dimensional non-unimodular Lorentzian Lie groups with Lie algebra $\g_7$ give solutions to the special Einstein-Weyl equation {\em (E-W)}, to the equation {\em (PS)} for a metric projective structure with a skew-symmetric Ricci tensor representative and to the vacuum near-horizon geometry equation {\em (VN-H)}. 
\end{cor}

\section{Generalized Ricci solitons on special $3D$ Lie groups}
\setcounter{equation}{0}

Let us consider a Lie algebra $\g$ (of dimension $\geq 2$) with the property that there exists a linear map $l: \g \to \mathbb R$, such that 
\begin{equation}\label{specLie}
[x, y] = l (x) y - l(y)x, \qquad x,y \in \g ,
\end{equation}
that is, the bracket product $[x, y]$ is always a linear
combination of $x$ and $y$. Left-invariant Riemannian metrics on a Lie algebra described by \eqref{specLie} were considered in \cite{M}, showing that they have constant sectional curvature $K=-||l||^2$ (thus, negative, unless $\g$ is abelian). Left-invariant Lorentzian metrics on Lie algebra \eqref{specLie} were investigated in \cite{No}. Again, they are of constant sectional curvature $K$, but $K$ can attain any real value. 

For a three-dimensional Lie algebra $\g$, with respect to any basis $\{e_1,e_2,e_3\}$,  from \eqref{specLie} we get
\begin{equation}\label{specRie}
[e_1,e_2] = B e_1-A e_2, \quad [e_1,e_3] = C e_1-A e_3, \quad [e_2,e_3] = C e_2-B e_3, 
\end{equation}
for three real constants $A=l(e_1),B=l(e_2),C=l(e_3)$. When we consider a left-invariant Riemannian metric $g$ on $\g$, we take 
$\{e_1,e_2,e_3\}$ as an orthonormal basis. Proceeding as in Section~2, we can then describe the Ricci tensor and $\mathcal{L}_X g$ with respect to $\{e_1,e_2,e_3\}$, for any left-invariant vector $X=X_i e_i$. In this way, \eqref{GRS} yields
\begin{equation}\label{sys10}
\left\{
\begin{array}{l}
B X_2+C X_3 +\alpha X_1^2+2\beta (A^2+B^2+C^2)=\lambda, \\[4pt]
A X_1+C X_3 +\alpha X_2^2+2\beta (A^2+B^2+C^2)=\lambda, \\[4pt]
A X_1 +B X_2+\alpha X_3^2+2\beta (A^2+B^2+C^2)=\lambda, \\[4pt]
-B X_1-A X_2  +2\alpha X_1X_2 =0, \\[4pt]
-C X_1-A X_3+2\alpha X_1X_3 =0, \\[4pt]
-C X_2-B X_3+2\alpha X_2X_3 =0.
\end{array}
\right.
\end{equation}

By a similar argument, for a left-invariant Lorentzian metric $g$ on $\g$, we take a pseudo-orthonormal basis $\{e_1,e_2,e_3\}$, with $e_3$ time-like, and $\g$ is again described by \eqref{specRie}. We then proceed as in Section~3 and find that for a left-invariant vector $X=X_i e_i$, equation \eqref{GRS} now becomes
\begin{equation}\label{sys11}
\left\{
\begin{array}{l}
B X_2+C X_3 +\alpha X_1^2+2\beta (A^2+B^2-C^2)=\lambda, \\[4pt]
A X_1+C X_3 +\alpha X_2^2+2\beta (A^2+B^2-C^2)=\lambda, \\[4pt]
-A X_1 -B X_2+\alpha X_3^2-2\beta (A^2+B^2-C^2)=-\lambda, \\[4pt]
-B X_1-A X_2  +2\alpha X_1X_2 =0, \\[4pt]
-C X_1+A X_3-2\alpha X_1X_3 =0, \\[4pt]
-C X_2+B X_3-2\alpha X_2X_3 =0.
\end{array}
\right.
\end{equation}
Only some changes of sign occur between systems \eqref{sys10} and \eqref{sys11}, due to the different signatures of the metrics. However, these slight changes are responsible, once again, of differences concerning their solutions. Solving \eqref{sys10} and \eqref{sys11}, we prove the following.

\begin{theorem}\label{GRSspec}
Let $\g$ denote a three-dimensional Lie algebra described by \eqref{specLie}. 

\smallskip
(I) If $g$ is a left-invariant Riemannian metric on $\g$, then \eqref{specRie} holds with respect to an orthonormal basis $\{e_1,e_2,e_3\}$, and the nontrivial left-invariant generalized Ricci solitons on $\g$ are the following:

\smallskip
(1) $\lambda = 2\beta (A^2+B^2+C^2)$, $X=0$: the metric is Einstein.

\smallskip
(2) $A=B=C=0$, $\alpha=0$, $\lambda=0$, $X$ is arbitrary: abelian Lie algebra.

\smallskip
(3) $\lambda=(\frac{1}{\alpha}+2\beta)(A^2+B^2+C^2)$, $X_1=\frac{A}{\alpha}$, $X_2=\frac{B}{\alpha}$, $X_3=\frac{C}{\alpha}$.  

\smallskip
(II) If $g$ is a left-invariant Lorentzian metric on $\g$, then \eqref{specRie} holds with respect to a pseudo-orthonormal basis $\{e_1,e_2,e_3\}$, with $e_3$ time-like, and the nontrivial left-invariant generalized Ricci solitons on $\g$ are the following:

\smallskip
(1')  $\lambda = 2\beta (A^2+B^2-C^2)$, $X=0$: the metric is Einstein.

\smallskip
(2') $A=B=C=0$, $\alpha=0$, $\lambda=0$, $X$ is arbitrary: abelian Lie algebra.

\smallskip
(3') either $A=B=0$ or $C=0$, and $\lambda=(\frac{1}{\alpha}+2\beta)(A^2+B^2-C^2)$, $X_1=\frac{A}{\alpha}$, $X_2=\frac{B}{\alpha}$, $X_3=-\frac{C}{\alpha}$.

\end{theorem}

\smallskip\noindent
It may be observed that solution {\em (3)} of the Riemannian case corresponds to solution {\em (3')} of the Lorentzian one, which holds either  for a space-like vector $X= \frac{A}{\alpha} e_1+ \frac{B}{\alpha} e_2$, or for a time-like vector $X_3=-\frac{C}{\alpha} e_3$.
  
It is easily seen that solutions {\em (3)} and {\em (3')} are compatible with (E-W) and (VN-H) and also with (PS), as $\frac{1}{\alpha}+2\beta=0$ when $\alpha=1 $ and $\beta=-\frac 12$. Therefore, the following result holds.

\begin{cor}
Three-dimensional  Lie groups, with special Lie algebra \eqref{specLie}, give solutions to the special Einstein-Weyl equation {\em (E-W)}, to the equation {\em (PS)} for a metric projective structure with a skew-symmetric Ricci tensor representative and to the vacuum near-horizon geometry equation {\em (VN-H)}. 
\end{cor}

\section{Final remarks}
\setcounter{equation}{0}

We end with some observations concerning the geometry of some of the solutions we found to the generalized Ricci soliton equation \eqref{GRS} for three-dimensional Lie groups. 

The most interesting case occurring in Theorem~\ref{GRSuniRie}, namely, solution {\em (4)}, requires the structure coefficients $A,B,C$ to satisfy $A=B\neq C$ (or any of their permutation). This condition occurs in the classification of naturally reductive Riemannian 
three-manifolds \cite{TV}. More precisely, for $C=0$ one gets a flat metric on $\widetilde{E}(2)$, while for $C\neq 0$, one finds the proper (that is, not locally symmetric) examples of naturally reductive spaces as left-invariant metrics on $SU(2)$, $\widetilde{SL}(2,\mathbb R)$ and $H_3$. Correspondingly, naturally reductive Lorentzian three-manifolds are found in the cases of Lie algebra $\g_3$ satisfying $A=B\neq C$ (with $C\neq 0$ in the proper case), and of Lie algebra $\g_4$ satisfying $A=B-\eta$ (with $A\neq 0$ in the proper case) (see \cite{CM0},\cite{CM}). Again, these conditions appear in the solutions to the generalized Ricci soliton equation for these Lorentzian Lie algebras, in Theorems~\ref{GRSg3} and \ref{GRSg4}, respectively.

If a Riemannian manifold $(M^n,g)$ admits a parallel line field $\mathcal D$, then $(M^n,g)$ splits locally into the direct product of a line and an $(n-1)$-dimensional manifold. The same property is true for a pseudo-Riemannian manifold admitting a parallel line field $\mathcal D$, either space-like or time-like. However, in pseudo-Riemannian settings, a different phenomenon can occur: it may exist a parallel degenerate line field $\mathcal D$, locally spanned by a light-like vector field $U$ satisfying $\nabla U=\omega \otimes U$.  Pseudo-Riemannian manifolds admitting a parallel vector field are a special case of the so called {\em Walker manifolds}, and are responsible for several behaviours which do not have any Riemannian counterpart. We may refer to \cite{CGV} for the investigation of three-dimensional Walker manifolds.
 
With regard to Lie algebra $\g_1$, when $B=0$ (that is, for Lorentzian Lie group $E(1,1)$), vector $X=X_2(e_2+e_3)$ is light-like and $\nabla X$ is parallel to $X$, for any real constant $X_2 \neq 0$. Therefore, this is a Walker manifold, and a vector of the form $X=X_2(e_2+e_3)$  occurs in solution {\em (3)} of Theorem~\ref{GRSg1} for the generalized Ricci soliton equation. For non-unimodular Lie algebra $\g_7$, whenever $C=0$, the Lorentzian metric is again Walker (with both locally symmetric and not locally symmetric examples occurring), as $\nabla X$ is parallel to $X$ when $X=X_2(e_2+e_3)$. And a vector of this form appears in solution {\em (3)} of Theorem~\ref{GRSg7} for the generalized Ricci soliton equation.


\begin{thebibliography}{9999}


\bibitem{BCGG}
M. Brozos-Vazquez, G. Calvaruso, E. Garcia-Rio and S. Gavino-Fernandez, \emph{Three-dimensional Lorentzian homogeneous Ricci solitons},
Israel J. Math. \textbf{188} (2012), 385--403.

\bibitem{CP}
D.M.J. Calderbank and H. Pedersen, \emph{Einstein–Weyl geometry}, In: LeBrun, C.,Wang, M. (Eds.) Surveys
in Differential Geometry, vol. VI: Essays on Einstein Manifolds, Suppl. to Journal of Differential
Geometry.

\bibitem{C}
G. Calvaruso, \emph{Homogeneous structures on three-dimensional Lo\-rentz\-ian manifolds}, J. Geom. Phys. {\bf 57} (2007), 1279--1291. Addendum: J. Geom. Phys. {\bf 58} (2008), 291--292.

\bibitem{C2}
G. Calvaruso, \emph{Einstein-like metrics on three-dimensional
homogeneous Lorentzian manifolds}, Geom. Dedicata {\bf 127}
(2007), 99--119.

\bibitem{CM0}
G. Calvaruso and R.A. Marinosci, \emph{Homogeneous geodesics of three-dimensional unimodular Lorentzian Lie groups}, Mediterr. J. Math. {\bf 3} (2006), 467--481.

\bibitem{CM}
G. Calvaruso and R.A. Marinosci, \emph{Homogeneous geodesics of non-unimodular Lorentz\-ian
Lie groups and naturally reductive Lorentzian spaces
in dimension three}, Adv. Geom. {\bf 8} (2008), 473--489.

\bibitem{Cao}
H.-D. Cao, \emph{Recent progress on Ricci solitons}, in Recent Advances in Geometric Analysis,
Advanced Lectures in Mathematics (ALM) 11, Int. Press, Somerville, MA, 2010, pp. 1–38.

\bibitem{Cerbo}
L.F. di Cerbo, \emph{Generic properties of homogeneous Ricci solitons},  Adv. Geom.  {\bf 14}  (2014), 225--237.

\bibitem{CGV}
 M. Chaichi, E. Garc\'{i}a-R\'{i}o and M.E. V\'{a}zquez-Abal, \emph{Three-dimensional Lorentz manifolds admitting a parallel null vector field}, {J. Phys. A: Math. Gen.} {\bf 38} (2005), 841--850.

\bibitem{CRT}
P.T. Chrusciel, H.S. Reall and P. Tod, \emph{On non-existence of static vacuum black holes with degenerate
components of the event horizon}, Class. Quantum Gravity {\bf 23} (2006), 549--554.

 \bibitem{CoPa}
L.A. Cordero and P.E. Parker, \emph{Left-invariant Lorentzian metrics on $3$-dimensional Lie groups}, Rend. Mat., Serie VII {\bf 17} (1997), 129--155.

\bibitem{HL}
K.Y. Ha and J.B. Lee, \emph{Left-invariant metrics and curvatures on simply connected three-dimensional Lie groups}, Math. Nachr. {\bf 282} (2009),  868--898.

\bibitem{M}
J. Milnor, \emph{Curvature of left invariant metrics on Lie groups}, Adv. Math. {\bf 21} (1976), 293--329.

\bibitem{No}
K. Nomizu, \emph{Left invariant Lorentz metrics on Lie groups}, Osaka J. Math. {\bf 16} (1979), 143--150.

\bibitem{NR}
P. Nurowski and M. Randall, \emph{Generalized Ricci solitons},
J. Geom. Anal., to appear. DOI:10.1007/s12220-015-9592-8.

\bibitem{R}
S. Rahmani, \emph{M\'{e}triques de Lorentz sur les groupes de Lie unimodulaires de dimension trois}, J. Geom. Phys. {\bf 9} (1992), 295--302.

\bibitem{Ra}
M. Randall, \emph{Local obstructions to projective surfaces admitting skew-symmetric Ricci tensor}, J. Geom.
Phys. {\bf 76} (2014), 192--199.

\bibitem{S}
K. Sekigawa, \emph{On some three-dimensional curvature homogeneous spaces}, Tensor N.S. {\bf 31} (1977), 87--97.

\bibitem{TV}
F. Tricerri and L. Vanhecke, \emph{Homogeneous structures on Riemannian manifolds}, London Mathematical Society Lecture Note Series, 
{\bf 83}, Cambridge Univ. Press (1983).
\end{thebibliography}
\end{document}